\documentclass{article}
\usepackage[english,russian,ukrainian]{babel}
\usepackage{latexsym,amssymb}
\usepackage[cp1251]{inputenc}
\usepackage{amsfonts,amssymb}
\usepackage{euscript}
\newcounter{theorem} %[section]
\newcounter{lemma} %[section]

\def\Re{\mathop\mathrm{Re}}\,
\def\Im{\mathop\mathrm{Im}}\,

\begin{document}

\vspace*{7mm}

\Large

\addtocounter{page}{-1}
\thispagestyle{empty}

\noindent УДК 517.548\\
\noindent\textbf{В.~С. Шпаківський, Т.~С. Кузьменко}
%\vspace{3mm}

\noindent\textbf{Про один клас кватерніонних відображень}

\vspace{7mm}

Рассмотрен новый класс кватернионных отображений, имеющих
связь с пространственными уравнениями в частных производных.
 Получено описание
всех отображений из этого класса, с помощью четырех аналитических
функций комплексной переменной.

\vspace{15mm}

\noindent {\bf  V.~S.~Shpakivskyi, T.~S.~Kuzmenko}\\
\noindent {\bf On one class of quaternionic mappings}

\vspace{7mm}

We consider a new class of quaternionic mappings,
associated with the spatial partial differential equations.
 We describe all mappings from this class using four analytic
  functions of the complex variable.

\Large
\newpage

\textbf{1. Вступ.}
Кватерніонний аналіз вже давно сформувався і активно розвивається як
окремий напрямок в математиці, завдяки його численним застосуванням в різних галузях
науки, переважно в математичній фізиці та в диференціальних рівняннях (див.,
наприклад, \cite{Gurlebeck,Kravchenko}). Реалізація такого зв'язку полягає
у введенні спеціальних класів кватерніонних ''диференційовних''\, функцій, компоненти
яких задовольняють певні системи диференціальних рівнянь типу системи Коші -- Рімана.
%В межах такого підходу, найбільш цікавим є розгляд систем диференціальних рівнянь
%в яких кількість змінних і рівнянь змінюються від 3 до 4.

Так початком кватерніонного аналізу у просторі $\mathbb{R}^3$ була робота Г.~Моісіла
 і Н.~Теодореско \cite{Moisil}, у якій вперше запропоновано тривимірний аналог
  системи рівнянь Коші -- Рімана. Вони ввели поняття \textit{голоморфного вектора},
    як кватерніоннозначної вектор-функції, компоненти якої неперервно диференційовні
  і задовольняють згадану вище систему, що дістала назву системи Моісіла -- Теодореско.
   В тій же роботі \cite{Moisil}
  автори довели аналог теореми Морера та аналоги інтегральної
  теореми та інтегральної формули Коші. Започатковані в \cite{Moisil}
  дослідження були продовжені в роботі \cite{Bitsadze}, де введено поняття інтеграла
  типу Коші та досліджено існування його граничних значень, а
  також знайдено його застосування до систем сингулярних
  інтегральних рівнянь.

Р.~Фютер \cite{Fueter} побудував чотиривимірне узагальнення системи Моісіла -- Теодореско
 та для введених ним \textit{регулярних}
функцій довів аналоги класичних результатів комплексного аналізу.

Згадані дослідження були узагальнені в роботі
\cite{Sudbery} і разом із застосуваннями у деяких моделях математичної фізики,
 відображені також в монографії \cite{Kravchenko}.
 Слід також відмітити, що так звані \textit{$\alpha$-голоморфні} функції
 $f$, які є об'єктом
 дослідження роботи \cite{Kravchenko}, задовольняють тривимірне рівняння Гельмгольца
 $$(\Delta_3+\alpha)f:=\frac{\partial^2 f}{\partial x^2}+\frac{\partial^2 f}
{\partial y^2}+
\frac{\partial^2 f}{\partial z^2}+\alpha f=0,
 $$
 де $\alpha$ --- кватерніон.

Останні дослідження у цьому напрямку (див., наприклад,
\cite{Gerus},\cite{Gerus-Shap},\cite{Schneider})
 полягають в різного роду узагальненнях
результатів роботи \cite{Kravchenko}.

Іншим, порівняно новим, напрямком кватерніонного аналізу в $\mathbb{R}^3$
і $\mathbb{R}^4$ є так
званий модифікований кватерніонний аналіз, започаткований
Г.~Льойтвілером на початку 90-х років (див., наприклад, [7 -- 9]).
У конструкції Г.~Льойтвілера в $\mathbb{R}^3$ перші дві компоненти, введених ним
\textit{гіперголоморфних} функцій
$f=u(x,y,z)+iv(x,y,z)+jw(x,y,z)$ \, (де $i,j$ --- базисні кватерніонні одиниці),
 задовольняють рівняння Лапласа-Бельтрамі
$$z\Delta_3u-\frac{\partial u}{\partial z}=0,
$$
 а третя компонента $w$ --- рівняння
$$z^2\Delta_3w-z\frac{\partial w}{\partial z}+w=0.
$$

В роботі \cite{Leutwiler-92} отримано розклад гіперголоморфної функції
в ряд по деякій системі кватерніонних поліномів.

На відміну від робіт \cite{Moisil,Fueter,Sudbery,Kravchenko},
в підході Г.~Льойтвілера гіперголорфною є
степенева функція, а частинні похідні гіперголорфної функції знову гіперголорфні.
В той же час між описаними вище напрямками існує певний зв'язок (див.
\cite{Eriksson-Bique}).

Ще однією сучасною теорією в кватерніонному аналізі є теорія так званих
\textit{$s$-регулярних} функцій, які введені Г.~Джентілі та Д.~Струппою
в роботі \cite{Gentili-Struppa-2006} в результаті розвитку ідеї К.~Кулліна \cite{Cullen}.

Ідея полягає в наступному. Кожен кватерніон $x:=x_0+\Im\,x$ \, при\,
 $x\neq x_0$
 можна подати у вигляді ''комплексного числа''\, з новою уявною одиницею $I$, а саме:
$x=x_0+I\,|\Im\,x|$,  де \, $I:=\frac{\Im\,x}{|\Im\,x|}$, а
 $|\cdot|$ --- модуль кватерніона. Очевидно, що $I^2=-1$.
  У такому ж вигляді можна подати й
 кватерніоннозначну функцію:\,
 $f(x)=U(x_0,|\Im\,x|)+I\,V(x_0,|\Im\,x|)$. Тоді функція $f$
 називається \textit{$s$-регулярною}  (див. \cite{Gentili-Struppa-2006}), якщо  ''комплекснозначна'' функція
$f=U+IV$ є голоморфною функцією ''комплексної'' змінної $x=x_0+I\,|\Im\,x|$.
 Очевидно, що $s$-регулярними є
 всі кватерніонні поліноми.
У наш час теорія $s$-регулярних функцій
  продовжує стрімко розвиватися
  (див. монографії \cite{Colombo-Sab-Struppa-2011,Gentili}).

В цій роботі розглядається спеціальний клас відображень в алгебрі
комплексних кватерніонів, який не охоплюється згаданими вище теоріями.
 Зазначимо, що комутативна
алгебра бікомплексних чисел (або комутативних кватерніонів Сегре \cite{Segre})
є підалгеброю  алгебри комплексних кватерніонів.
%Алгебра Сегре може бути розглянута як двовимірна комутативна напівпроста алгебра
%над полем комплексних чисел $\mathbb{C}$.
В цій підалгебрі виділимо
тривимірний дійсний підпростір і  розглянемо відображення, які визначені в
області цього підпростору і приймають значення у всій алгебрі
комплексних кватерніонів. Такі відображення, які є неперервними і
диференційовними за Гато, назвемо \textit{$G$-моногенними}. Вони і є
 основним об'єктом дослідження.

Встановлено, що $G$-моногенними є не лише кватерніонні поліноми, а й
кватерніонні степеневі ряди. Більше того, в роботі встановлено конструктивний опис
усіх $G$-моногенних відображень за допомогою чотирьох аналітичних
функцій комплексної змінної. Як наслідок, похідна Гато $G$-моногенного
відображення в свою чергу є $G$-моногенним відображенням.
 Крім того, досліджено зв'язок $G$-моногенних відображень
з просторовими рівняннями з частинними похідними. Зокрема, наведено застосування
моногенних відображень до побудови розв'язків тривимірного рівняння Лапласа.

\textbf{2. Алгебра комплексних кватерніонів.}
Нехай $\mathbb{H(C)}$ --- алгебра кватерніонів над полем комплексних чисел
$\mathbb{C}$, базис якої складається з одиниці алгебри $1$ і елементів $I,J,K,$
для яких виконуються
наступні правила
множення:
$$
I^2=J^2=K^2=-1,\,IJ=-JI=K,\,JK=-KJ=I,\,KI=-IK=J.
$$

Розглянемо в алгебрі $\mathbb{H(C)}$ інший базис $\{e_1,e_2,e_3,e_4\}$, розклад
 елементів якого в базисі $\{1,I,J,K\}$ має вигляд
$$e_1=\frac{1}{2}(1+iI), \quad e_2=\frac{1}{2}(1-iI), \quad
e_3=\frac{1}{2}(iJ-K), \quad e_4=\frac{1}{2}(iJ+K),$$
де $i$ --- уявна комплексна одиниця. Таблиця множення в новому базисі
набуває вигляду

\begin{equation}\label{tabl}
\begin{tabular}{c||c|c|c|c|}
$\cdot$ & $e_1$ & $e_2$ & $e_3$ & $e_4$\\
\hline
\hline
$e_1$ & $e_1$ & $0$ & $e_3$ & $0$\\
\hline
$e_2$ & $0$ & $e_2$ & $0$ & $e_4$\\
\hline
$e_3$ & $0$ & $e_3$ & $0$ & $e_1$\\
\hline
$e_4$ & $e_4$ & $0$ & $e_2$ & $0$\\
\hline
\end{tabular}\,\,.
\end{equation}\vskip2mm

Норма кватерніона $a=\sum\limits_{k=1}^4a_ke_k,\,a_k\in\mathbb{C}$
визначається рівністю
\begin{equation}\label{norma-1}
\|a\|:=\sqrt{\sum\limits_{k=1}^4|a_k|^2}\,,
\end{equation}
а
одиниця алгебри $\mathbb{H(C)}$ в цьому базисі є сумою ідемпотентів: $1=e_1+e_2$.
Очевидно також, що комутативна підалгебра
 з базисом $\{e_1,e_2\}$ є згаданою вище алгеброю бікомплексних чисел або
 алгеброю комутативних кватерніонів Сегре \cite{Segre}.

Нагадаємо (див., наприклад, \cite[c. 64]{van_der_Varden}), що підмножина
$\mathcal{I}\subset\mathbb{H(C)}$
називається \textit{лівим} (або \textit{правим}) \textit{ідеалом}, якщо з
 умови $x\in\mathcal{I}$
випливає $yx\in\mathcal{I}$\, (або $xy\in\mathcal{I}$) для довільного
$y\in\mathbb{H(C)}$.
Тепер відмітимо, що алгебра $\mathbb{H(C)}$ містить два праві максимальні ідеали
$$\mathcal{I}_1:=\{\lambda_2e_2+\lambda_4e_4:\lambda_2,\lambda_4\in\mathbb{C}\},\qquad
\mathcal{I}_2:=\{\lambda_1e_1+\lambda_3e_3:\lambda_1,\lambda_3\in\mathbb{C}\}$$
і два ліві максимальні ідеали
$$\mathcal{\widehat{I}}_1:=\{\lambda_2e_2+\lambda_3e_3:\lambda_2,\lambda_3\in
\mathbb{C}\},\qquad
\mathcal{\widehat{I}}_2:=\{\lambda_1e_1+\lambda_4e_4:\lambda_1,\lambda_4\in\mathbb{C}
\}.$$
Очевидно, що радикал алгебри складається лише з нульового елемента, тобто алгебра
$\mathbb{H(C)}$ напівпроста.

Наслідком очевидних рівностей
$$\mathcal{I}_1\cap\mathcal{I}_2=\widehat{\mathcal{I}}_1\cap\widehat{\mathcal{I}}_2=
0,\qquad \mathcal{I}_1\cup\mathcal{I}_2=\widehat{\mathcal{I}}_1\cup
\widehat{\mathcal{I}}_2=\mathbb{H(C)}
$$
є розклад в пряму суму:
$$\mathbb{H(C)}=\mathcal{I}_1\oplus\mathcal{I}_2=
\widehat{\mathcal{I}}_1\oplus\widehat{\mathcal{I}}_2.
$$

Введемо в розгляд лінійні функціонали $f_1:\mathbb{H(C)\rightarrow\mathbb{C}}$ та $f_2:\mathbb{H(C)\rightarrow\mathbb{C}},$ покладаючи
$$f_1(e_1)=f_1(e_3)=1, \qquad f_1(e_2)=f_1(e_4)=0,$$
$$f_2(e_2)=f_2(e_4)=1, \qquad f_2(e_1)=f_2(e_3)=0,$$
при цьому очевидно $f_1(\mathcal{I}_1)=f_2(\mathcal{I}_2)=0$.

Визначимо також лінійні функціонали
$\widehat{f}_1:\mathbb{H(C)\rightarrow\mathbb{C}}$ та $\widehat{f}_2:\mathbb{H(C)\rightarrow\mathbb{C}}$ рівностями
$$\widehat{f}_1(e_1)=\widehat{f}_1(e_4)=1, \qquad \widehat{f}_1(e_2)=
\widehat{f}_1(e_3)=0,$$
$$\widehat{f}_2(e_2)=\widehat{f}_2(e_3)=1, \qquad \widehat{f}_2(e_1)=
\widehat{f}_2(e_4)=0,$$
для яких очевидно $\widehat{f}_1(\widehat{\mathcal{I}}_1)=
\widehat{f}_2(\widehat{\mathcal{I}}_2)=0$.

\textbf{3. Ліво-$G$-моногенні та право-$G$-моногенні відображення.}
Нехай
\begin{equation}\label{i-bazis}
i_1=1, \qquad i_2=a_1e_1+a_2e_2, \qquad i_3=b_1e_1+b_2e_2
\end{equation}
при $a_k,b_k\in\mathbb{C},\,\,
k=1,2$ --- трійка лінійно незалежних векторів над полем дійсних чисел $\mathbb{R}$.
Це означає, що рівність $$\alpha_1i_1+\alpha_2i_2+\alpha_3i_3=0, \qquad
\alpha_1,\alpha_2,\alpha_3\in\mathbb{R}$$ виконується тоді і тільки тоді,
коли  $\alpha_1=\alpha_2=\alpha_3=0$.

Виділимо в алгебрі $\mathbb{H(C)}$ лінійну оболонку $E_3:=\{\zeta=xi_1+yi_2+zi_3:x,y,z\in\mathbb{R}\}$ над полем дійсних чисел
$\mathbb{R},$ породжену векторами $i_1,i_2,i_3.$ Області $\Omega$
тривимірного простору $\mathbb{R}^3$ поставимо у відповідність область $\Omega_\zeta:=\{\zeta=xi_1+yi_2+zi_3:(x,y,z)\in\Omega\}$ в $E_3$.

Введемо позначення $$\xi_1:=f_1(\zeta)=\widehat{f}_1(\zeta)=x+ya_1+zb_1,$$
$$\xi_2:=
f_2(\zeta)=\widehat{f}_2(\zeta)=x+ya_2+zb_2.$$
%За визначенням  функціоналів, числа
% $f_1,f_2,\widehat{f}_1,\widehat{f}_2$,\,
% $\xi_1$ і $\xi_2$ --- комплексні.
Тепер елемент $\zeta\in E_3$ можна подати у
 вигляді $\zeta=\xi_1e_1+\xi_2e_2$, і згідно визначення (\ref{norma-1})
\begin{equation}\label{norma}
\|\zeta\|=\sqrt{|\xi_1|^2+|\xi_2|^2}.
\end{equation}

Відмітимо, що в подальшому  істотним є припущення:
 $f_1(E_3)=f_2(E_3)=\mathbb{C}$.
Очевидно, що воно має місце тоді і тільки тоді, коли
хоча б одне з чисел у кожній з пар $(a_1,b_1)$, $(a_2,b_2)$ належить
 $\mathbb{C}\setminus\mathbb{R}$.

Скажемо, що деякий функціонал $f:\mathbb{H(C)}\rightarrow\mathbb{C}$
\textit{право--муль\-ти\-плі\-ка\-тив\-ний} (або \textit{ліво--мультиплікативний}),
 якщо
для довільних $x\in\mathbb{H(C)}$ і $y\in E_3$ справедлива рівність
$f(yx)=f(y)f(x)$\, \big(або $f(xy)=f(x)f(y)$\big).
\vskip 2mm
\textbf{Лема 1.} \emph{Функціонали $f_1:\mathbb{H(C)\rightarrow\mathbb{C}}$ та $f_2:\mathbb{H(C)}\rightarrow\mathbb{C}$ неперервні і
право--мультиплікативні, а функціонали $\widehat{f}_1:\mathbb{H(C)
\rightarrow\mathbb{C}}$ та $\widehat{f}_2:\mathbb{H(C)\rightarrow\mathbb{C}}$
 неперервні і
ліво--мультиплікативні.}
\vskip 1mm

\textbf{Доведення.} Відповідна мультиплікативність всіх
 функціоналів встановлюється безпосередньою перевіркою, а неперервність ---
 випливає з їх обмеженості. А саме, якщо $a=\sum\limits_{k=1}^4a_ke_k\in
 \mathbb{H(C)}$,
  то, наприклад для $f_1$, маємо
 $$\frac{|f_1(a)|}{\|a\|}\leq\frac{|a_1|+|a_3|}{\sqrt{|a_1|^2+|a_2|^2+
 |a_3|^2+|a_4|^2}}\leq2.$$
  Аналогічно доводиться неперервність інших функціоналів.
  Лему доведено.

Неперервне відображення $\Phi:\Omega_\zeta\rightarrow\mathbb{H(C)}$ (або $\widehat{\Phi}:\Omega_\zeta\rightarrow\mathbb{H(C)}$) називається
\emph{право-$G$-моногенним}
\big(або \emph{ліво-$G$-моногенним}\big) в області
$\Omega_\zeta\subset E_3$, якщо $\Phi$ \big(або $\widehat{\Phi}$\big)
диференційовне за Гато у кожній точці цієї області, тобто якщо для
кожного $\zeta\in\Omega_\zeta$ існує елемент
$\Phi'(\zeta)$ \big(або $\widehat{\Phi}'(\zeta)$\big) алгебри
$\mathbb{H(C)}$ такий,
що виконується рівність
\begin{equation}\label{ozn-l-monog}
\lim\limits_{\varepsilon\rightarrow 0+0}\Big(\Phi(\zeta+\varepsilon h)-\Phi(\zeta)\Big)\varepsilon^{-1}= h\Phi'(\zeta)\quad\forall\,h\in E_3
\end{equation}

\begin{equation}\label{ozn-r-monog}
\Biggr(\mbox{або}\,\, \lim\limits_{\varepsilon\rightarrow 0+0}
\left(\widehat{\Phi}(\zeta+\varepsilon h)-\widehat{\Phi}(\zeta)\right)
\varepsilon^{-1}= \widehat{\Phi}'(\zeta)h\quad\forall\,h\in E_3\Biggr).
\end{equation}
При цьому $\Phi'(\zeta)$ назвемо \emph{правою похідною Гато} в точці $\zeta$\,, а
  $\widehat{\Phi}'(\zeta)$ ---
 \emph{лівою похідною Гато} в точці $\zeta$\,.

\vskip 1mm
\textbf{Теорема 1.}\emph{ Відображення $\Phi:\Omega_\zeta\rightarrow\mathbb{H(C)}$
 вигляду
\begin{equation}\label{Phi}
\Phi(\zeta)=\sum_{k=1}^{4}{U_k(x,y,z)e_k}, \quad x,y,z\in \mathbb{R},
\end{equation}
де $U_k:\Omega\rightarrow\mathbb{C}$ --- диференційовні функції в
 області $\Omega$, є
право-$G$-моногенним  або ліво-$G$-моногенним в області
$\Omega_\zeta\subset E_3$ тоді і
тільки тоді, коли виконуються відповідно умови:}

$$\frac{\partial U_1}{\partial y}=a_1\frac{\partial U_1}{\partial x},\quad
\frac{\partial U_2}{\partial y}=a_2\frac{\partial U_2}{\partial x},\quad
 \frac{\partial U_3}{\partial y}=a_2\frac{\partial U_3}{\partial x},\quad
\frac{\partial U_4}{\partial y}=a_1\frac{\partial U_4}{\partial x},$$\vspace{-10mm}
\begin{equation}\label{teor-r-K-R}
\end{equation}
 $$\frac{\partial U_1}{\partial z}=b_1\frac{\partial U_1}{\partial x},\quad
\frac{\partial U_2}{\partial z}=b_2\frac{\partial U_2}{\partial x},
\quad \frac{\partial U_3}{\partial z}=b_2\frac{\partial U_3}{\partial x},
\quad \frac{\partial U_4}{\partial z}=b_1\frac{\partial U_4}{\partial x}.$$%\\
\emph{або}
$$\frac{\partial U_1}{\partial y}=a_1\frac{\partial U_1}{\partial x},\quad
\frac{\partial U_2}{\partial y}=a_2\frac{\partial U_2}{\partial x},\quad
 \frac{\partial U_3}{\partial y}=a_1\frac{\partial U_3}{\partial x},\quad
\frac{\partial U_4}{\partial y}=a_2\frac{\partial U_4}{\partial x},$$\vspace{-10mm}
\begin{equation}\label{teor-l-K-R}
\end{equation}
 $$\frac{\partial U_1}{\partial z}=b_1\frac{\partial U_1}{\partial x},\quad
\frac{\partial U_2}{\partial z}=b_2\frac{\partial U_2}{\partial x},
\quad \frac{\partial U_3}{\partial z}=b_1\frac{\partial U_3}{\partial x},
\quad \frac{\partial U_4}{\partial z}=b_2\frac{\partial U_4}{\partial x}.$$
\vskip 3mm

\textbf{Доведення.} \emph{Необхідність.} Якщо відображення (\ref{Phi})
право-$G$-моногенне в області $\Omega_\zeta$,
 то при $h=i_1$ рівність (\ref{ozn-l-monog}) набуває вигляду
$$\Phi'(\zeta)=\sum_{k=1}^{4}{\frac{\partial U_k(x,y,z)}{\partial x}\,e_k}, \quad \zeta=xi_1+yi_2+zi_3\in\Omega_\zeta.$$
Тепер, покладаючи в рівності (\ref{ozn-l-monog}) спочатку $h=i_2,$ а потім
 $h=i_3,$ та з урахуванням правил множення для базисних елементів, отримаємо умови (\ref{teor-r-K-R}) для компонент право-$G$-моногенного відображення (\ref{Phi}).

\emph{Достатність.} Нехай $\zeta=xi_1+yi_2+zi_3\in\Omega_\zeta,
\,h:=h_1i_1+h_2i_2+h_3i_3,$ де $h_1,h_2,h_3\in\mathbb{R}$ і
додатне число $\varepsilon$ таке,
 що $\zeta+\varepsilon h\in\Omega_\zeta.$ Враховуючи умови (\ref{teor-r-K-R}),
  маємо
$$\Big(\Phi(\zeta+\varepsilon h)-\Phi(\zeta)\Big)\varepsilon^{-1}-h
\sum_{k=1}^{4}{\frac{\partial U_k(x,y,z)}{\partial x}\,e_k}=$$
$$=\varepsilon^{-1}\sum_{k=1}^{4}{\Big(U_k(x+\varepsilon h_1,y+
\varepsilon h_2,z+\varepsilon h_3)-U_k(x,y,z)\Big)e_k-}$$
$$-\left(\frac{\partial U_1}{\partial x}h_1+a_1\frac{\partial U_1}
{\partial x}h_2+b_1\frac{\partial U_1}{\partial x}h_3\right)e_1-
\left(\frac{\partial U_2}{\partial x}h_1+a_2\frac{\partial U_2}
{\partial x}h_2+b_2\frac{\partial U_2}{\partial x}h_3\right)e_2-$$
$$-\left(\frac{\partial U_3}{\partial x}h_1+a_1\frac{\partial U_3}
{\partial x}h_2+b_1\frac{\partial U_3}{\partial x}h_3\right)e_3-
\left(\frac{\partial U_4}{\partial x}h_1+a_2\frac{\partial U_4}
{\partial x}h_2+b_2\frac{\partial U_4}{\partial x}h_3\right)e_4=$$
$$=\varepsilon^{-1}\sum_{k=1}^{4}{\Biggr(U_k(x+\varepsilon h_1,y+
\varepsilon h_2,z+\varepsilon h_3)-U_k(x,y,z)-}$$
\begin{equation}\label{peretv}
-\frac{\partial U_k(x,y,z)}{\partial x}\,\varepsilon h_1-
\frac{\partial U_k(x,y,z)}{\partial y}\,\varepsilon h_2-
\frac{\partial U_k(x,y,z)}{\partial z}\,\varepsilon h_3\Biggr)e_k.
\end{equation}

Внаслідок диференційовності функцій $U_k$ в області $\Omega$
справедливі співвідношення
$$U_k(x+\varepsilon h_1,y+\varepsilon h_2,z+\varepsilon h_3)-
U_k(x,y,z)-\frac{\partial U_k(x,y,z)}{\partial x}\varepsilon h_1-$$
$$-\frac{\partial U_k(x,y,z)}{\partial y}\varepsilon h_2-\frac
{\partial U_k(x,y,z)}{\partial z}\varepsilon h_3=o(\varepsilon),
\quad \varepsilon\rightarrow 0,\,\, k=\overline{1,4}.$$

Тому, перейшовши до границі в рівності (\ref{peretv}) при
 $\varepsilon\rightarrow 0,$ отримаємо рівність (\ref{ozn-l-monog}).
  Аналогічно доводиться випадок для ліво-$G$-моногенного відображення.
   Теорему доведено.

Відмітимо, що умови (\ref{teor-r-K-R}) і (\ref{teor-l-K-R})  є аналогами
умов Коші -- Рімана
і у згорнутому вигляді можуть бути записані так:
\begin{equation}\label{umova-l-K-R}
\frac{\partial \Phi}{\partial y}=i_2\,\frac{\partial \Phi}{\partial x},
\qquad \frac{\partial \Phi}{\partial z}=i_3\,\frac{\partial \Phi}{\partial x}
\end{equation}
для право-$G$-моногенного відображення, і
\begin{equation}\label{umova-r-K-R}
\frac{\partial \widehat{\Phi}}{\partial y}=\frac{\partial
\widehat{\Phi}}{\partial x}\,i_2, \qquad \frac{\partial
\widehat{\Phi}}{\partial z}=\frac{\partial \widehat{\Phi}}{\partial x}\,i_3
\end{equation}
для ліво-$G$-моногенного відображення.

Розглянемо приклади право- і ліво-$G$-моногенних відображень.
Враховуючи подання елемента $\zeta$ у вигляді $\zeta=\xi_1e_1+\xi_2e_2$ і таблицю
множення алгебри $\mathbb{H(C)}$, маємо \,
$\zeta^n=\xi_1^n\,e_1+\xi_2^n\,e_2$. Шляхом перевірки умов
(\ref{umova-l-K-R}), (\ref{umova-r-K-R}), легко переконатися в тому,
що відображення $\Phi(\zeta)=\zeta^n$ є одночасно право- і ліво-$G$-моногенним
у всьому просторі $E_3$.
Аналогічно перевіряється, що відображення
\begin{equation}\label{222}
\Phi(\zeta)=\sum\limits_
{k=0}^n\zeta^k\,c_k,\quad c_k\in\mathbb{H(C)}
\end{equation}
є право-$G$-моногенним в $E_3$, а відображення
$$\widehat{\Phi}(\zeta)=\sum\limits_
{k=0}^nc_k\,\zeta^k,\quad c_k\in\mathbb{H(C)}$$
--- ліво-$G$-моногенним в $E_3$.

\textbf{4. Конструктивний опис право-$G$-моногенних і
 ліво-$G$-моногенних відображень.}

\textbf{Лема 2.} \textit{Розклад резольвенти має вигляд}
\begin{equation}\label{lema-resolventa}
(t-\zeta)^{-1}=\frac{1}{t-\xi_1}\,e_1+\frac{1}{t-\xi_2}\,e_2,\quad
\forall\,\,t\in \mathbb{C}:\,t\neq \xi_1, \, t\neq \xi_2.
\end{equation}

\textbf{Доведення.} Встановимо при яких $t\in\mathbb{C}$ в алгебрі
$\mathbb{H(C)}$ існує елемент $(t-\zeta)^{-1}$ і знайдемо коефіцієнти $A_k$ його
розкладу за базисом: $$(t-\zeta)^{-1}=\sum\limits_{k=1}^4A_k\,e_k.$$

Враховуючи подання (\ref{i-bazis}) елементів $i_1,i_2,i_3$ за базисом
$\{e_1,e_2,e_3,e_4\}$ і таблицю множення алгебри $\mathbb{H(C)}$, маємо
$$1=(t-\zeta)(t-\zeta)^{-1}=\Big((t-\xi_1)e_1+(t-\xi_2)e_2\Big)
\sum\limits_{k=1}^4A_k\,e_k=
$$
$$=(t-\xi_1)A_1e_1+(t-\xi_1)A_3e_3+(t-\xi_2)A_2e_2+(t-\xi_2)A_4e_4=e_1+e_2.
$$
Тепер прирівнюючи коефіцієнти при відповідних базисних одиницях,
отримуємо розклад (\ref{lema-resolventa}). Лему доведено.

Із рівності (\ref{lema-resolventa})
випливає, що точки $(x,y,z)\in\mathbb{R}^3,$ які відповідають необоротним
елементам $\zeta=xi_1+yi_2+zi_3\in E_3$, лежать на прямих
$$L_1: x+y\Re\,a_1+z\Re\,b_1=0,
\qquad y\Im\,a_1+z\Im\,b_1=0,$$
$$L_2: x+y\Re\,a_2+z\Re\,b_2=0,
\qquad y\Im\,a_2+z\Im\,b_2=0$$
в просторі $\mathbb{R}^3$.

Область $\Omega\in\mathbb{R}^3$ називають \emph{опуклою в напрямку прямої $L,$}
якщо вона містить кожен відрізок, який паралельний прямій $L$ і з'єднує дві
 точки цієї області.

\vspace{1mm}
\textbf{Лема 3.} \emph{Нехай область $\Omega\in\mathbb{R}^3$ опукла в
 напрямку прямих $L_1$ і $L_2$,\,\,\, $f_1(E_3)=f_2(E_3)=\mathbb{C}$, а
 відображення $\Phi:\Omega_\zeta\rightarrow\mathbb{H(C)}$
 право-$G$-моногенне в області $\Omega_\zeta$. Якщо точки
  $\zeta_1,\zeta_2\in\Omega_\zeta$ такі, що $\zeta_1-\zeta_2\in
  \{\zeta=xi_1+yi_2+zi_3:(x,y,z)\in L_1\},$ то}
\begin{equation}\label{osn-lema-1}
\Phi(\zeta_1)-\Phi(\zeta_2)\in\mathcal{I}_1.
\end{equation}
\emph{Якщо ж точки $\zeta_1,\zeta_2\in\Omega_\zeta$ такі, що $\zeta_1-\zeta_2\in\{\zeta=xi_1+yi_2+zi_3:(x,y,z)\in L_2\},$ то}
\begin{equation}\label{osn-lema-2}
\Phi(\zeta_1)-\Phi(\zeta_2)\in\mathcal{I}_2.
\end{equation}
\vskip 1mm

Співвідношення (\ref{osn-lema-1}) доводиться за схемою доведення леми 1 роботи \cite{Plaksa_Shpakivskyi}, в якому замість прямої $L$ необхідно використати пряму
 $L_1,$ а замість функціонала $f$ потрібно використовувати функціонал $f_1$.
 Аналогічно доводиться співвідношення (\ref{osn-lema-2}) з заміною $L_1$ і
 $f_1$ відповідно на $L_2$ і $f_2$.
 При цьому використовується лема 1 цієї роботи.

Повністю аналогічно доводиться і наступне твердження.

\vskip 1mm
\textbf{Лема 4.} \emph{Нехай область $\Omega\in\mathbb{R}^3$
опукла в напрямку прямих $L_1$ і $L_2$,\,\,\, $f_1(E_3)=
f_2(E_3)=\mathbb{C}$, а відображення $\widehat{\Phi}:\Omega_\zeta\rightarrow
\mathbb{H(C)}$
ліво-$G$-моногенне в області $\Omega_\zeta$. Якщо
 точки $\zeta_1,\zeta_2\in\Omega_\zeta$ такі, що $\zeta_1-\zeta_2\in\{\zeta=xi_1+yi_2+zi_3:(x,y,z)\in L_1\},$ то}
$$\widehat{\Phi}(\zeta_1)-\widehat{\Phi}(\zeta_2)\in\widehat{\mathcal{I}}_1.$$
\emph{Якщо ж точки $\zeta_1,\zeta_2\in\Omega_\zeta$ такі, що $\zeta_1-\zeta_2\in\{\zeta=xi_1+yi_2+zi_3:(x,y,z)\in L_2\},$ то}
$$\widehat{\Phi}(\zeta_1)-\widehat{\Phi}(\zeta_2)\in\widehat{\mathcal{I}}_2.$$
\vskip 1mm

\vskip 1mm
\textbf{Теорема 2.} \emph{Кожне право-$G$-моногенне в
області $\Omega_\zeta$ відображення $\Phi:\Omega_\zeta\rightarrow
\mathbb{H(C)}$ подається у вигляді}
$$\Phi(\zeta)=\Phi_{10}(\zeta)+\Phi_{20}(\zeta),$$
\emph{де $\Phi_{10}:\Omega_\zeta\rightarrow\mathcal{I}_1,\,\,\,\Phi_{20}:\Omega_
\zeta\rightarrow\mathcal{I}_2$ --- деякі право-$G$-моногенні в області
$\Omega_\zeta$ відображення зі значеннями в правих максимальних ідеалах $\mathcal{I}_1,\,\mathcal{I}_2,$ а кожне ліво-$G$-моногенне в
 області $\Omega_\zeta$ відображення подається у вигляді}
\begin{equation}\label{1111}
\widehat{\Phi}(\zeta)=\widehat{\Phi}_{10}(\zeta)+\widehat{\Phi}_{20}(\zeta),
\end{equation}
\emph{де $\widehat{\Phi}_{10}:\Omega_\zeta\rightarrow\widehat{\mathcal{I}}_1,
\,\,\,\widehat{\Phi}_{20}:\Omega_\zeta\rightarrow\widehat{\mathcal{I}}_2$ --- деякі ліво-$G$-моногенні в області $\Omega_\zeta$ відображення зі
 значеннями в лівих максимальних ідеалах $\widehat{\mathcal{I}}_1,\,
 \widehat{\mathcal{I}}_2$.}
\vskip 1mm

\textbf{Доведення.} Із розкладу одиниці
 $1=e_1+e_2$ випливає, що довільне відображення
 $\Phi:\Omega_\zeta\rightarrow\mathbb{H(C)}$ подається у вигляді
  $$\Phi=e_1\Phi+e_2\Phi$$
  і при цьому $e_1\Phi\in\mathcal{I}_2$, а $e_2\Phi\in\mathcal{I}_1$.

Введемо позначення $\Phi_{10}:=e_2\Phi$, \, $\Phi_{20}:=e_1\Phi$. Покажемо,
 що відображення $\Phi_{10},\,\Phi_{20}$  право-$G$-моногенні
  в області $\Omega_\zeta$. Для цього рівність (\ref{ozn-l-monog}) помножимо
  зліва на $e_1$:
\begin{equation}\label{teor-Phi}
\lim\limits_{\varepsilon\rightarrow 0+0}e_1\Big(\Phi(\zeta+\varepsilon h)-\Phi(\zeta)\Big)\varepsilon^{-1}= e_1h\Phi'(\zeta)\quad\forall\,h\in E_3.
\end{equation}
Оскільки елементи $e_1$ та $h$ належать комутативній підалгебрі з
базисом $\{e_1,e_2\},$ то $e_1h=he_1$, і тому з рівності
(\ref{teor-Phi}) випливає рівність
$$
\lim\limits_{\varepsilon\rightarrow 0+0}\Big(e_1\Phi(\zeta+\varepsilon h)-e_1\Phi(\zeta)\Big)\varepsilon^{-1}= he_1\Phi'(\zeta),
$$
яка і доводить, що відображення $\Phi_{20}$
право-$G$-моногенне в області $\Omega_\zeta$. Аналогічно доводиться
права-$G$-моногенність
відображення $\Phi_{10}$.

Аналогічно доводиться представлення  (\ref{1111}).
 Теорему доведено.

В наступній теоремі описано усі право- та ліво-$G$-моногенні
відображення зі значеннями відповідно в ідеалах $\mathcal{I}_1$
  та $\widehat{\mathcal{I}}_1$ за допомогою аналітичних
  функцій відповідної комплексної змінної.

  Для формулювання результату введемо позначення
   $$ D_1:=f_1(\Omega_\zeta)\subset\mathbb{C} ,
  \quad  D_2:=f_2(\Omega_\zeta)\subset \mathbb{C}.$$

 \vskip 1mm
\textbf{Теорема 3.} \emph{Нехай область $\Omega$ опукла в напрямку прямої
$L_2$ і $f_1(E_3)=f_2(E_3)=\mathbb{C}$.
Тоді кожне право-$G$-моногенне відображення $\Phi_{10}:\Omega_
\zeta\rightarrow\mathcal{I}_1$  подається у вигляді
\begin{equation}\label{Phi10-l}
\Phi_{10}(\zeta)=F_{12}(\xi_2)e_2+F_{14}(\xi_2)e_4 \qquad \forall\,\zeta=xi_1+yi_2+zi_3\in\Omega_\zeta,
\end{equation}
а ліво-$G$-моногенне відображення $\widehat{\Phi}_{10}:\Omega_\zeta\rightarrow\mathcal{\widehat{I}}_1$
подається у такому вигляді}
\begin{equation}\label{Phi10-r}
\widehat{\Phi}_{10}(\zeta)=F_{11}(\xi_2)e_2+F_{13}(\xi_2)e_3 \qquad \forall\,\zeta=xi_1+yi_2+zi_3\in\Omega_\zeta,
\end{equation}
\emph{де $F_{12},F_{14},F_{11},F_{13}$ --- деякі аналітичні в області
$D_2$ функції змінної $\xi_2:=x+ya_2+zb_2$.}
\vskip 1mm

\textbf{Доведення.} Оскільки $\Phi_{10}$ приймає значення в ідеалі
$\mathcal{I}_1,$ то справедлива рівність
\begin{equation}\label{Phi10-I1}
\Phi_{10}(\zeta)=V_2(x,y,z)e_2+V_4(x,y,z)e_4,
\end{equation}
де $V_2:\Omega\rightarrow\mathbb{C}$ і $V_4:\Omega\rightarrow\mathbb{C}.$

Для відображення $\Phi_{10}$ виконуються умови правої-$G$-
моногенності (\ref{umova-l-K-R}) при $\Phi=\Phi_{10},$ з
яких після підстановки в них виразів (\ref{i-bazis}),
(\ref{Phi10-I1}), з урахуванням однозначності розкладу
елементів алгебри $\mathbb{H(C)}$ за базисом $\{e_1,e_2,e_3,e_4\}$,
 отримаємо систему для знаходження функцій $V_2,V_4:$
\begin{equation}\label{syst-Phi10}
\begin{array}{l}
\displaystyle \frac{\partial V_2}{\partial y}=a_2\frac{\partial V_2}{\partial x}\,,\\[5mm]
\displaystyle\frac{\partial V_4}{\partial y}=a_2\frac{\partial V_4}{\partial x}\,,\\[5mm]
\displaystyle\frac{\partial V_2}{\partial z}=b_2\frac{\partial V_2}{\partial x}\,,\\[5mm]
\displaystyle\frac{\partial V_4}{\partial z}=b_2\frac{\partial V_4}{\partial x}.
\end{array}
\end{equation}

З першого і третього рівняння системи (\ref{syst-Phi10}) знайдемо функцію $V_2$.
Для цього виділимо дійсну і уявну частину змінної $\xi_2$:
$$\xi_2=(x+y\Re\,a_2+z\Re\,b_2)+i(y\Im\,a_2+
z\Im\,b_2):=\tau_2+i\eta_2$$
і відмітимо, що наслідком вказаних рівнянь є рівності
\begin{equation}\label{211}
\frac{\partial V_2}{\partial\eta_2}\,\Im\,a_2
=i\frac{\partial V_2}{\partial\tau_2}\,\Im\,a_2,\quad
\frac{\partial V_2}{\partial\eta_2}\,\Im\,b_2
=i\frac{\partial V_2}{\partial\tau_2}\,\Im\,b_2.
\end{equation}

Оскільки $f_1(E_3)=f_2(E_3)=\mathbb{C}$ випливає, що хоча б одне з чисел
$\Im\,a_2$ або $\Im\,b_2$ відмінне від нуля, то з
(\ref{211}) отримуємо рівність
$$\frac{\partial V_2}{\partial\eta_2}=i\frac{\partial V_2}{\partial\tau_2}.$$

Тепер так, як і при доведенні теореми 2 з
\cite{Plaksa_Shpakivskyi}, з використанням леми 3 і теореми 6 з
\cite{Tolstov} доводиться рівність
$V_2(x_1,y_1,z_1)=V_2(x_2,y_2,z_2)$ для точок $(x_1,y_1,z_1),
(x_2,y_2,z_2)\in\Omega$ таких, що відрізок, який з'єднує ці
 точки, паралельний прямій $L_2.$ Звідси випливає, що функція $V_2$ вигляду $V_2(x,y,z):=F_{12}(\xi_2),$ де $F_{12}$ --- довільна аналітична
 функція в області $D_2,$ є загальним розв'язком системи, яка
 складається з першого і третього рівнянь системи (\ref{syst-Phi10}).

Тепер з другого і четвертого рівнянь системи (\ref{syst-Phi10})
аналогічно встановлюємо, що функція $V_4$ має
 вигляд $V_4(x,y,z):=F_{14}(\xi_2),$ де $F_{14}$ --- довільна аналітична
в області $D_2$ функція.

Повністю аналогічно доводиться рівність (\ref{Phi10-r}). Теорему доведено.

В наступній теоремі, яка доводиться повність аналогічно до
 теореми 3, описано усі право- та ліво-$G$-моногенні
відображення зі значеннями відповідно в ідеалах $\mathcal{I}_2$
 та $\mathcal{\widehat{I}}_2$ алгебри $\mathbb{H(C)}$ за допомогою
  аналітичних функцій відповідної комплексної змінної.

\vskip 1mm
\textbf{Теорема 4.} \emph{Нехай область $\Omega$ опукла в напрямку прямої
$L_1$ і $f_1(E_3)=f_2(E_3)=\mathbb{C}$.
Тоді кожне право-$G$-моногенне відображення $\Phi_{20}:\Omega_
\zeta\rightarrow\mathcal{I}_2$  подається у вигляді
\begin{equation}\label{Phi20-l}
\Phi_{20}(\zeta)=F_{21}(\xi_1)e_1+F_{23}(\xi_1)e_3 \qquad \forall\,\zeta=xi_1+yi_2+zi_3\in\Omega_\zeta,
\end{equation}
а ліво-$G$-моногенне відображення $\widehat{\Phi}_{20}:\Omega_\zeta\rightarrow\mathcal{\widehat{I}}_2$
подається у такому вигляді}
\begin{equation}\label{Phi20-r}
\widehat{\Phi}_{20}(\zeta)=F_{22}(\xi_2)e_1+F_{24}(\xi_1)e_4 \qquad \forall\,\zeta=xi_1+yi_2+zi_3\in\Omega_\zeta,
\end{equation}
\emph{де $F_{21},F_{23},F_{22},F_{24}$ --- деякі аналітичні в області
$D_1$ функції змінної $\xi_1:=x+ya_1+zb_1$.}
\vskip 1mm

З урахуванням теореми 2 та рівностей (\ref{Phi10-l}), (\ref{Phi20-l})
 доведено наступне твердження.

\vskip 1mm
\textbf{Теорема 5.} \emph{Нехай область $\Omega$ опукла в
напрямку прямих $L_1$ і $L_2,$ а $f_1(E_3)=f_2(E_3)=\mathbb{C}$.
 Тоді кожне право-$G$-моногенне відображення
 $\Phi:\Omega_\zeta\rightarrow\mathbb{H(C)}$ подається у вигляді}
\begin{equation}\label{Phi-l-rozklad}
\Phi(\zeta)=F_1(\xi_1)e_1+F_2(\xi_2)e_2+F_3(\xi_1)e_3+F_4(\xi_2)e_4
\end{equation}
$$\forall\,\zeta=xi_1+yi_2+zi_3\in\Omega_\zeta,$$
\emph{де $F_1$ і $F_3$ --- деякі аналітичні в області $D_1$
 функції змінної $\xi_1:=x+ya_1+zb_1$, а $F_2$ і $F_4$ ---
 деякі аналітичні в області $D_2$ функції змінної $\xi_2:=x+ya_2+zb_2$.}
\vskip 1mm

Тепер очевидно, що розглянуте вище відображення (\ref{222}) буде
право-$G$-моногенним в $E_3$, оскільки для нього функції
$F_1,F_2,F_3,F_4$ будуть поліномами. Але тепер можна сказати й
більше. А саме,  право-$G$-моногенним у відповідній області
відображенням  буде
 не тільки поліном вигляду (\ref{222}), а й ряд вигляду
\begin{equation}\label{222-1}
\Phi(\zeta)=\sum\limits_
{k=0}^\infty\zeta^k\,c_k,\qquad c_k\in\mathbb{H(C)},
\end{equation}
для якого комплексні степеневі ряди, що виступають в ролі
 аналітичних фунцій $F_1,F_2,F_3,F_4$, є  збіжними.

Тепер, в силу рівностей (\ref{Phi10-r}) і (\ref{Phi20-r}),
 справедливе твердження для ліво-$G$-моногенного відображення.

\vskip 1mm
\textbf{Теорема 6.} \emph{Нехай область $\Omega$ опукла в
напрямку прямих $L_1$ і $L_2,$ а $f_1(E_3)=f_2(E_3)
=\mathbb{C}$. Тоді кожне ліво-$G$-моногенне відображення $\widehat{\Phi}:\Omega_\zeta\rightarrow\mathbb{H(C)}$ подається у вигляді}
\begin{equation}\label{Phi-r-rozklad}
\widehat{\Phi}(\zeta)=F_1(\xi_1)e_1+F_2(\xi_2)e_2+F_3(\xi_2)e_3+F_4(\xi_1)e_4
\end{equation}
$$\forall\,\zeta=xi_1+yi_2+zi_3\in\Omega_\zeta,$$
\emph{де $F_1$ і $F_4$ --- деякі аналітичні в області $D_1$
функції змінної $\xi_1:=x+ya_1+zb_1$, а $F_2$ і $F_3$ --- деякі
аналітичні в області $D_2$ функції змінної $\xi_2:=x+ya_2+zb_2.$}
\vskip 1mm

Аналогічно до (\ref{222-1}),  відображення
\begin{equation}\label{222-2}
\widehat{\Phi}(\zeta)=\sum\limits_
{k=0}^\infty c_k\,\zeta^k,\qquad c_k\in\mathbb{H(C)},
\end{equation}
є ліво-$g$-моногенним.

Очевидно, що формула (\ref{Phi-l-rozklad}) дає можливість
побудувати усі право-$G$-моногенні відображення $\Phi:\Omega_\zeta
\rightarrow\mathbb{H(C)}$, а формула (\ref{Phi-r-rozklad})
дає можливість побудувати усі ліво-$G$-моногенні
відображення $\widehat{\Phi}:\Omega_\zeta\rightarrow\mathbb{H(C)}$
за допомогою чотирьох аналітичних функцій відповідної комплексної змінної.
Відмітимо, що в роботі \cite{Flaut-Shpakiv-Hyperhol} за допомогою аналітичних функцій
 комплексної змінної
будуються так звані $A_t$\textit{-гіперголоморфні} функції в довільній алгебрі
Келі-Діксона $A_t$ над полем $\mathbb{R}$.

Порівнюючи праві частини рівностей (\ref{Phi-l-rozklad}) і
(\ref{Phi-r-rozklad}), приходимо до висновку, що відображення
$\Psi(\zeta)$ буде одночасно право-
 і ліво-$G$-моногенним тоді і тільки тоді, коли воно має вигляд
 $\Psi(\zeta)=F_1(\xi_1)e_1+F_2(\xi_2)e_2+c_3e_3+c_4e_4$, де $c_3,c_4\in\mathbb{C}$.
  Тепер очевидно, що відображеня
 $\Psi(\zeta)=\zeta^n=\xi_1^n\,e_1+\xi_2^n\,e_2$ одночасно
 право-  і ліво-$G$-моногенне в $E_3$.

Прямим наслідком представлень (\ref{Phi-l-rozklad}),
(\ref{Phi-r-rozklad}) є той факт, що множини усіх право- і
ліво-$G$-моногенних відображень зі
 значення в алгебрі $\mathbb{H(C)}$ утворюють
функціональні алгебри в областях їх визначення. Тобто добуток
двох, наприклад,
 право-$G$-моногенних відображень знову є право-$G$-моногенним відображенням.

З урахуванням розкладу (\ref{lema-resolventa}) і правил множення (\ref{tabl}),
 одержуємо наступні інтегральні представлення
  право- і ліво-$G$-моногенних відображень
$$\Phi(\zeta)=\frac{1}{2\pi i}\int\limits_{\Gamma_1}(t-\zeta)^{-1}
\Big(F_1(t)e_1+F_3(t)e_3\Big)dt+
\frac{1}{2\pi i}\int\limits_{\Gamma_2}(t-\zeta)^{-1}
\Big(F_2(t)e_2+F_4(t)e_4\Big)dt,$$

$$\widehat{\Phi}(\zeta)=\frac{1}{2\pi i}\int\limits_{\Gamma_1}
\Big(F_1(t)e_1+F_4(t)e_4\Big)
(t-\zeta)^{-1}
dt+
\frac{1}{2\pi i}\int\limits_{\Gamma_2}
\Big(F_2(t)e_2+F_3(t)e_3\Big)(t-\zeta)^{-1}
dt,$$
де замкнені жорданові спрямлювані криві $\Gamma_k$ лежать у
відповідних областях $D_k$, охоплюють відповідні точки
$\xi_k$ і не містять точок $\xi_n$,\, $k,n=1,2$ при $k\neq n$.

Відмітимо також, що похідні Гато право-$G$-моногенного відображення
$\Phi(\zeta)$ і ліво-$G$-моногенного відображення $\widehat{\Phi}
(\zeta)$ виражаються відповідно формулами
\begin{equation}\label{Phi-r----}
\Phi'(\zeta)=F_1'(\xi_1)e_1+F_2'(\xi_2)e_2+F_3'(\xi_1)e_3+F_4'
(\xi_2)e_4
\end{equation}
і
$$\widehat{\Phi}'(\zeta)=F_1'(\xi_1)e_1+F_2'(\xi_2)e_2+F_3'(\xi_2
)e_3+F_4'(\xi_1)e_4.$$

Наступне твердження випливає безпосередньо з рівності (\ref{Phi-l-rozklad}) і (\ref{Phi-r-rozklad}), праві частини яких є відповідно право- і
 ліво-$G$-моногенним відображенням в області $\Pi_\zeta:=\{\zeta\in E_3:f_1(\zeta)\in D_1,\,f_2(\zeta)\in D_2\}.$

\vskip 1mm
\textbf{Теорема 7.} \emph{Нехай область $\Omega$ опукла в
 напрямку прямих $L_1$ і $L_2$,  $f_1(E_3)= f_2(E_3)=\mathbb{C},$ і відображення $\Phi:\Omega_\zeta\rightarrow\mathbb{H(C)}$ право-$G$-моногенне, а $\widehat{\Phi}:\Omega_\zeta\rightarrow\mathbb{H(C)}$ ---
 ліво-$G$-моногенне в області $\Omega_\zeta$. Тоді $\Phi$ і $\widehat{\Phi}$
 продовжуються до відображеннь, які є відповідно право- і
  ліво-$G$-моногенним в області $\Pi_\zeta.$}
\vskip 1mm

Теорема 7 дає можливість легко знайти область право-$G$-моногенності
відображення (\ref{222-1}) і ліво-$G$-моногенності відображення (\ref{222-2}).

Принциповим наслідком рівностей (\ref{Phi-l-rozklad}) та
(\ref{Phi-r-rozklad}) є наступне твердження, справедиве для
право- і  ліво-$G$-моногенних
 відображень в довільній області $\Omega_\zeta.$

\vskip 1mm
\textbf{Теорема 8.} \emph{Нехай область $\Omega$ опукла в
 напрямку прямих $L_1$ і $L_2$,  $f_1(E_3)= f_2(E_3)=\mathbb{C},$ і
  відображення $\Phi:\Omega_\zeta\rightarrow\mathbb{H(C)}$ право-$G$-моногенне, а $\widehat{\Phi}:\Omega_\zeta\rightarrow\mathbb{H(C)}$ ---
 ліво-$G$-моногенне в області $\Omega_\zeta$. Тоді
 похідні Гато усіх порядків відображеннь $\Phi$ і $\widehat{\Phi}$
  є відповідно право- і ліво-$G$-моногенними відображеннями
  в області $\Omega_\zeta$.}
\vskip 1mm

\textbf{Доведення.} Оскільки куля $\mho$ (яка повністю міститься
в області $\Omega$) з центром в довільній точці $(x_0,y_0,z_0)
\in\Omega$ є опуклою множиною в напрямку прямих $L_1$ і $L_2$,
 то в околі $\mho_\zeta:=\{\zeta=xi_1+yi_2+
zi_3:(x,y,z)\in\mho\}$ точки $\zeta_0=x_0i_1+y_0i_2+z_0i_3$ справедливі
  рівності (\ref{Phi-l-rozklad}) і (\ref{Phi-r----}). Але при цьому
 компоненти розкладу (\ref{Phi-r----}) є аналітичними функціями відповідних
  комплексних змінних,
 тобто вираз для $\Phi'(\zeta)$ має вигляд рівності (\ref{Phi-l-rozklad}), а
 це і означає праву-$G$-моногенність відображення  $\Phi'(\zeta)$.
 Аналогічно доводиться для похідної Гато довільного порядку
 і для ліво-$G$-моногенного відображення $\widehat{\Phi}(\zeta)$.
 Теорему доведено.

\textbf{5. Зв'язок право- і ліво-$G$-моногенних відображень з
рівняннями в частинних похідних.} Розглянемо лінійне
диференціальне рівняння в частинних
 похідних із сталими коефіцієнтами:
\begin{equation}\label{dif-riv}
\mathcal{L}_nU(x,y,z):=\sum\limits_{\alpha+\beta+\gamma=n}
C_{\alpha,\beta,\gamma}\frac{\partial^nU}{\partial x^\alpha\partial
y^\beta\partial z^\gamma}=0,\qquad C_{\alpha,\beta,\gamma}\in\mathbb{R}.
\end{equation}

Якщо відображення $\Phi(\zeta)$ є $n$ разів право-диференційовним
за Гато, а відображення $\widehat{\Phi}(\zeta)$ є $n$ разів
ліво-диференційовним за Гато,
 то наслідком рівностей (\ref{ozn-l-monog}) і (\ref{ozn-r-monog})
 є відповідно рівності
$$\frac{\partial^{\alpha+\beta+\gamma}\Phi}{\partial
x^\alpha\partial y^\beta\partial z^\gamma}=i_1^\alpha\, i_2^\beta\, i_3^\gamma\,\Phi^{(\alpha+\beta+\gamma)}(\zeta)=i_2^\beta\,
 i_3^\gamma\,\Phi^{(n)}(\zeta)$$
і
$$\frac{\partial^{\alpha+\beta+\gamma}\widehat{\Phi}}
{\partial x^\alpha\partial y^\beta\partial z^\gamma}=\widehat{\Phi}^{(\alpha+\beta+\gamma)}(\zeta)\,i_1^\alpha\, i_2^\beta\, i_3^\gamma=\widehat{\Phi}^{(n)}(\zeta)\,i_2^\beta\, i_3^\gamma.$$
Тому внаслідок рівності
\begin{equation}\label{dif-riv-1}
\mathcal{L}_n\Phi(\zeta)=\sum
\limits_{\alpha+\beta+\gamma=n}C_{\alpha,\beta,\gamma}\,i_2^\beta\,
 i_3^\gamma\,\,\Phi^{(n)}(\zeta)
\end{equation}
кожне $n$ разів право-диференційовне за Гато відображення $\Phi$ при
виконанні умов $\Phi^{(n)}(\zeta)\neq0$ і
\begin{equation}\label{sum}
\sum\limits_{\alpha+\beta+\gamma=n}C_{\alpha,\beta,\gamma}\,i_2^\beta\, i_3^\gamma=0
\end{equation}
задовольняє рівняння $\mathcal{L}_n\Phi(\zeta)=0$. Аналогічно
внаслідок рівності
\begin{equation}\label{dif-riv-2}
\mathcal{L}_n\widehat{\Phi}(\zeta)=\widehat{\Phi}^{(n)}(\zeta)
\sum\limits_{\alpha+\beta+\gamma=n}
C_{\alpha,\beta,\gamma}\,i_2^\beta\, i_3^\gamma
\end{equation}
кожне $n$ разів ліво-диференційовне за Гато відображення
$\widehat{\Phi}$ при виконанні умов
$\widehat{\Phi}^{(n)}(\zeta)\neq0$ і (\ref{sum})
задовольняє рівняння $\mathcal{L}_n\widehat{\Phi}(\zeta)=0.$ Відповідно,
усі дійснозначні компоненти розкладу відображень $\Phi$
і $\widehat{\Phi}$ за базисом $\{e_1,e_2,e_3,e_4,ie_1,ie_2,ie_3,ie_4\}$ є
розв'язками рівняння (\ref{dif-riv}).

Таким чином,
задача про побудову розв'язків рівняння (\ref{dif-riv})
у вигляді компонент право- або ліво-диференційовних за Гато
відображень зводиться до відшукання в алгебрі $\mathbb{H(C)}$
трійки лінійно незалежних над полем $\mathbb{R}$ векторів (\ref{i-bazis}), які
задовольняють характеристичне рівняння (\ref{sum}).

Відмітимо, що якщо обидва функціонали $f_1,\,f_2$ приймають
значення в $\mathbb{C}$, то згідно з теоремою 8 кожне  право- і
ліво-$G$-моногенне відображення задовольняє рівність
(\ref{dif-riv-1}).

Очевидно, що співвідношення
\begin{equation}\label{funk}
f_1(E_3)=f_2(E_3)=\mathbb{C}
\end{equation}
має місце тоді і тільки тоді, коли
хоча б одне з чисел у кожній з пар $(a_1,b_1)$, $(a_2,b_2)$ належить
 $\mathbb{C}\setminus\mathbb{R}$.
Якщо рівняння (\ref{dif-riv}) має особливий вигляд,
то можна вказати достатні умови для виконання співвідношень
(\ref{funk}). Для цього введемо позначення
\begin{equation}\label{P(a,b)}
P(a,b):=\sum\limits_{\alpha+\beta+\gamma=n}C_{\alpha,\beta,\gamma}\,
a^\beta\, b^\gamma.
\end{equation}

\vskip 1mm
\textbf{Теорема 9.} \emph{Нехай в алгебрі $\mathbb{H(C)}$ існує
 трійка лінійно незалежних над полем $\mathbb{R}$ векторів вигляду
  \em(\ref{i-bazis}),\em\, які задовольняють рівність \em(\ref{sum}).
  \em\, Тоді якщо $P(a,b)\neq0$ при всіх дійсних значеннях
  $a,b,$ то виконуються співвідношення} (\ref{funk}).
\vskip 1mm

\textbf{Доведення.} Використовуючи таблицю множення алгебри,
 маємо рівності
$$i_2^\beta=a_1^\beta e_1+a_2^\beta e_2, \qquad i_3^\gamma=
b_1^\gamma e_1+b_2^\gamma e_2.$$
Тепер рівність (\ref{sum}) набуває вигляду
$$
\sum\limits_{\alpha+\beta+\gamma=n}C_{\alpha,\beta,\gamma}
\left(a_1^\beta\, b_1^\gamma \,e_1+a_2^\beta\, b_2^\gamma\, e_2\right)=0.
$$
або в рівносильній формі
\begin{equation}\label{teor-9}
\sum\limits_{\alpha+\beta+\gamma=n}C_{\alpha,\beta,\gamma}\,
a_k^\beta\, b_k^\gamma =0,\qquad k=1,2.
\end{equation}

Оскільки розв'язок системи (\ref{teor-9}) існує (за умовою теореми) і $P(a,b)\neq0$
при всіх дійсних $a,b$, то рівності
(\ref{teor-9}) можуть виконуватися лише тоді, коли хоча б одне з
чисел у кожній з пар $(a_1,b_1)$, $(a_2,b_2)$ належить множині
$\mathbb{C}\backslash\mathbb{R}.$
  Теорему доведено.

Тепер зауважимо, що з умови теореми $P(a,b)\neq0$ випливає,
що завжди $C_{n,0,0}\neq0$, оскільки
в іншому випадку при $a=b=0$ було б $P(a,b)=0$. А також оскільки функція
$P(a,b)$ неперервна на $\mathbb{R}^2$, то умова $P(a,b)\neq0$ по суті
означає одне з двох  $P(a,b)>0$ або $P(a,b)<0$ при всіх $a,b\in\mathbb{R}$.

Очевидно також, що рівняння вигляду (\ref{dif-riv}) еліптичного типу
завжди задовольняє умову $P(a,b)\neq0$ при всіх $a,b\in\mathbb{R}$. Але в той же
час існують рівняння вигляду (\ref{dif-riv}) для яких $P(a,b)>0$ і
які не є еліптичними. Таким, наприклад, є рівняння
$$\frac{\partial^5 u}{\partial x^5}+
\frac{\partial^5 u}{\partial x\partial y^2\partial z^2}+
\frac{\partial^5 u}{\partial x\partial z^4}=0.
$$

\textbf{6. Приклад.} Покажемо зв'язок право- і ліво-$G$-моногенних
відображень з тривимірним рівнянням Лапласа:
\begin{equation}\label{Lapl}
\Delta_3U(x,y,z):=\frac{\partial^2 U}{\partial
x^2}+\frac{\partial^2 U} {\partial y^2}+ \frac{\partial^2
U}{\partial z^2}=0.
\end{equation}
Для рівняння (\ref{Lapl}) характеристичне
 рівняння (\ref{sum}), набуває вигляду
\begin{equation}\label{harm}
1+i_2^2+i_3^2=0.
\end{equation}

Подібно до \cite{Ketchum-28}, трійку лінійно назалежних над полем
$\mathbb{R}$ векторів $i_1=1,i_2,i_3$ назвемо \textit{гармонічною
трійкою}, якщо має місце рівність (\ref{harm}) і виконуються умови
$i_2^2\neq0$,  $i_3^2\neq0$.

Після підстановки рівностей (\ref{i-bazis}) в умови (\ref{harm})
приходимо до наступного твердження: \emph{гармонічними трійками в
алгебрі $\mathbb{H(C)}$ є вектори $1,i_2,i_3,$ розклад яких за
базисом $\{e_1,e_2,e_3,e_4\}$ має вигляд \em(\ref{i-bazis}) \em і
комплексні числа  $a_k, b_k$,\, $k=1,2$ задовольняють систему
рівнянь}
\begin{equation}\label{syst--1}
1+a_1^2+b_1^2=0,\qquad 1+a_2^2+b_2^2=0.
\end{equation}\vskip1mm

Систему (\ref{syst--1}) задовольняють, зокрема, вирази $a_1=i\sin
t$, $b_1=i\cos t$, $a_2=i\sin\tau$, $b_2=i\cos\tau$,\,
$t,\tau\in\mathbb{C}$\,, яким відповідають
\begin{equation}\label{1--1}
\xi_1=x+iy\sin t+iz\cos t,\quad \xi_2=x+iy\sin\tau+iz\cos\tau, \quad
t,\tau\in\mathbb{C}.
\end{equation}

Оскільки для рівняння Лапласа $P(a,b)=1+a^2+b^2>0$, то умови теореми 9
виконуються, а значить кожне право- і ліво-$G$-моногенне відображення
задовольняє рівняння (\ref{Lapl}). Представлення (\ref{Phi-l-rozklad}) і
 (\ref{Phi-r-rozklad}), в яких $\xi_1,\xi_2$ визначені рівностями (\ref{1--1}),
 визначають моногенні відображення в $\mathbb{H(C)}$, пов'язані з рівнянням
 (\ref{Lapl}). Звідси випливає, що розв'язками рівняння
 (\ref{Lapl}) є дійсна і уявна частини
функції $U(x,y,z)=F(x+iy\sin t+iz\cos t)$, де $t\in\mathbb{C}$
  і $F$ --- довільна аналітична функція.
%  На завершення відмітимо схожість отриманого розв'язку з тим, що наведено в роботі \cite[c. 238]{Uitteker-Vatson}.

\renewcommand{\refname}{Література}

\end{document}